\documentclass[12pt,a4paper]{article}
\usepackage {amssymb,latexsym}
\usepackage {amsmath}
\usepackage {stmaryrd}

\newtheorem{theorem}{Theorem}[section]
\newtheorem{proposition}{Proposition}[section]
\newtheorem{remark}{Remark}[section]
\newtheorem{lemma}{Lemma}[section]
\newtheorem{corollary}{Corollary}[section]

\newtheorem{example}{Eexample}[section]

\newenvironment{proof}[1][Proof]{ \textbf{#1 : }}
{\rule{0.5em}{0.5em}}

\newtheorem{definition}{Definition}

\newcommand{\ben}{\begin{enumerate}}
\newcommand{\een}{\end{enumerate}}

\begin{document}

\author{Bouaazza ES SAADI}

\title{A certain class of standard subalgebras of affine Kac-Moody algebras}
\date{}
\maketitle

\begin{center}
Laboratoire de Math\'{e}matiques et Applications, \\
Universit\'{e} de Haute Alsace, Mulhouse, France \vspace{4mm}
\\Email :\\
 Bouaazza.Es-saadi@uha.fr
\end{center}



\begin{abstract}
The aim of this paper is to extend the theory of standard
subalgebras of finite dimensional simple Lie algebras to infinite
dimensional Lie algebras. We construct and characterize a class of
standard subalgebras of affine Kac-Moody algebra.
\end{abstract}

AMS Classification 2000 : 17B20, 17B30, 17B65, 17B67.



\section*{Introduction}
Let $\mathfrak{g}$ be a finite dimensional simple Lie algebra over
$\mathbb{C}$, of rank $p$ and let $\mathfrak{h}$ be a Cartan
subalgebra of $\mathfrak{g}$. We denote by $\Delta$ the roots system
of the couple $(\mathfrak{g}, \mathfrak{h})$ and by $\Pi
=\{\alpha_{1}, \cdots, \alpha_{p}\}$ the simple roots of $\Delta$.

Let $A$ be a generalized Cartan matrix of affine type of order $p+1$
associated to finite dimensional simple Lie algebra $\mathfrak{g}$.

Let $S(A)$ (resp. $S(^tA)$) be the Dynkin diagram associated to $A$
(resp. $^tA$) and $\mathfrak{B}=\{\alpha_{0}, \cdots, \alpha_{p}\}$
(resp. $\mathfrak{B}^{\vee}=\{\alpha_{0}^{\vee}, \cdots,
\alpha_{p}^{\vee}\}$) be the set of the vertex of diagram $S(A)$
(resp. $S(^tA)$) called roots basis (resp. coroots basis). Let
$a_{0}, a_{1}, \cdots, a_{p}$ be the numerical labels of $S(A)$ and
$a_{0}^{\vee}, \cdots, a_{p}^{\vee}$ be the numerical labels of
$S(^tA)$.

We denote by $\mathcal{L}:=\mathbb{C}[t,t^{-1}]$ the algebra of
Laurent polynomials in $t$. Recall that the residue of a Laurent
polynomial $P= \underset{i\in \mathbb{Z}}{\sum }c_{i}t^{i}$ (where
all but a finite number of $c_{i}$ are $0$) is defined by $Res P =
c_{-1}$.

Consider the infinite dimensional Lie algebra, called \textit{the
Loop algebra}
$$\mathcal{L}(\mathfrak{g}):=\mathcal{L} \otimes
\mathfrak{g}$$ We denote by $[.,.]$ the bracket on $\mathfrak{g}$.
The bracket on $\mathcal{L}(\mathfrak{g})$ is defined as follows,
for all $(x,y) \in \mathfrak{g} \times \mathfrak{g}$ :
$[t^{n}\otimes x, t^{m} \otimes y]_{\mathcal{L}} = t^{n+m} \otimes
[x,y]$.

Fix a nondegenerate invariant symmetric bilinear $\mathbb{C}$-valued
form $(./.)$ on $\mathfrak{g}$. We extend this form by linearity to
an $\mathcal{L}$-valued bilinear form $(./.)_{L}$ on
$\mathcal{L}(\mathfrak{g})$ which is defined by
$$(P \otimes x/Q \otimes y)_{L}=PQ(x/y)$$

Now, we define a $\mathbb{C}$-valued $2$-cocycle on the Lie algebra
$\mathcal{L}(\mathfrak{g})$ by
$$ \psi(a,b)=Res(\frac{da}{dt},b)_{t} \qquad \text{for all} \; (a,b)
\in \mathcal{L}(\mathfrak{g}) \times \mathcal{L}(\mathfrak{g})$$ We
denote by $\tilde{\mathcal{L}}(\mathfrak{g})$ the extension of the
Lie algebra $\mathcal{L}(\mathfrak{g})$ by $1$-dimensional center,
associated to the cocycle $\psi$. Namely,
$\tilde{\mathcal{L}}(\mathfrak{g}) = \mathcal{L}(\mathfrak{g})
\oplus \mathbb{C} \mathrm{K}$ where $\mathrm{K} =
\underset{i=0}{\overset{i=p}{\sum}}a^{\vee}_{i}\alpha^{\vee}_{i}$ is
a generator of the $1$-dimensional center.

Finally, the affine Kac-Moody algebra associated to matrix $A$ is
$$\mathfrak{g}(A) := \mathcal{L}(\mathfrak{g}) \oplus \mathbb{C}
\mathrm{K} \oplus \mathbb{C} \texttt{d}$$ where $\texttt{d}$ is the
derivation defined by : $\texttt{d}=t\frac{d}{dt}$ on
$\mathcal{L}(\mathfrak{g})$, $\texttt{d}(\mathrm{K})=0$. For more
details of this construction see \cite[chap. 7]{VK}.

The aim of this paper is to construct a certain class of standard
subalgebras of affine Kac-Moody algebras $\mathfrak{g}(A)$. This
uses standard subalgebras of the complex simple Lie algebra
$\mathfrak{g}$.

Let $\tau$ be a subalgebra of $\mathfrak{g}$. \textit{ If the
normalizer of $\tau$ is a parabolic subalgebra of $\mathfrak{g}$
then $\tau$ is called standard subalgebra of $\mathfrak{g}$}.
Several papers were devoted to the study of standard subalgebras for
finite dimensional case, see \cite{GU}, \cite{KH},
\cite{EKM},\cite{KOP}, \cite{OP}.

These subalgebras were characterized in \cite{KH} and \cite{EKM}
using a noncomparable roots with respect to the partial order
relation $" \leqslant "$ defined on the dual vector space
$\mathfrak{h}^{*}$ by $\omega_{1} \leqslant \omega_{2}$ if
$\omega_{2} - \omega_{1}$ is linear combination of simple roots of
$\Pi$ with non-negative coefficients.

Let $\mathfrak{R}$ be a set of positive roots two by two
noncomparable. Consider the set $\mathfrak{R}_{1}=\{\beta \in
\Delta_{+}\, : \, \omega \leqslant \beta \, \text{ for certain }
\omega \in \mathfrak{R}\}$. Then the subalgebra $\mathfrak{m}
=\underset{\beta \in \mathfrak{R}_{1}}{\sum} \mathfrak{g}_{\beta}$
is nilpotent standard subalgebra of $\mathfrak{g}$.

Setting $\rho(\mathfrak{m})$ be the normalizer of $\mathfrak{m}$.
Let $\mathfrak{r}$ be an ideal of the Levi reductive subalgebra
$\mathfrak{h} + \underset{\alpha \in
\Omega_{1}}{\sum}\mathfrak{g}_{\alpha}$, lying in the Levi reductive
subalgebra associated to parabolic subalgebra $\rho(\mathfrak{m})$
where $\Omega_{1}$ is the set of roots which are linear combination
of $\Pi \backslash (\Pi \cap \mathfrak{R})$.

Any standard subalgebra $\tau$ of $\mathfrak{g}$ has the form
$$\tau = \mathfrak{m} + \mathfrak{r}$$ and the normalizer
$\rho(\tau)=\rho(\mathfrak{m})$.

In this paper, theorem \ref{existencetheorem} states that there
exists standard subalgebra $\overline{\tau}$ of $\mathfrak{g}(A)$
associated to a given standard subalgebra $\tau$ of finite
dimensional simple Lie algebra $\mathfrak{g}$ and a certain vector
subspace $V_{\tau}=[\tau, \mathfrak{g}]$ of $\mathfrak{g}$
associated to $\tau$. The standard subalgebra $\overline{\tau}$ has
the form
$$\overline{\tau} =  t^{n} \otimes \tau + t^{n+1} \otimes
V_{\tau} + t^{n+2}\mathbb{C}[t] \otimes \mathfrak{g} + \mathbb{C}K
\quad \text{ for } n \in \mathbb{N}^{*}
$$
The subspace $V_{\tau}$ is called "appui" subspace of
$\mathfrak{g}$.

The standard subalgebra $\overline{\tau}$ is determined by the
finite dimensional standard subalgebra $\tau$ and the appui subspace
$V_{\tau}$. Explicit formulas of $V_{\tau}$ are established for
different cases of $\tau$.

In the case where $\tau$ is nilpotent standard subalgebra. Then,
$\tau$ is associated to subset $\mathfrak{R} \subset \Delta_{+}$ of
positive roots two by two noncomparable. The theorem
\ref{theorem3.1} states that the appui subspace $V_{\tau}$ has the
following form :
$$V_{\tau} = \underset{\alpha \in \Delta _{+}\backslash
<\mathfrak{R}_{3}>^{+}}{\sum }\mathfrak{g}_{\alpha } +
\underset{\alpha \in \mathfrak{R}_{1}\cup \mathfrak{R}_{2}}{\sum}
\mathbb{C}\mathbb{[}e_{\alpha },e_{-\alpha }]+\underset{\alpha \in
\mathfrak{R}_{2}}{\sum }\mathfrak{g}_{-\alpha }$$ where
$\mathfrak{R}_{1}$, $\mathfrak{R}_{2}$ and $\mathfrak{R}_{3}$ are
defined by the relation \ref{R1}, \ref{R2} and \ref{R3}.

In the second case where $\tau$ is standard subalgebra of nilpotent
radical. Then, $\tau = \mathfrak{m} + \mathfrak{r}_{0}$ where
$\mathfrak{m}$ is nilpotent part of $\tau$ and $\mathfrak{r}_{0}$ is
semisimple part of $\tau$. Therefor, $\mathfrak{r}_{0}$ is
constructed from a certain common connected subsystem $\psi$ of $\Pi
\backslash (\Pi \cap \mathfrak{R})$ and $\Pi \backslash
(\underset{\omega \in \mathfrak{R}}{\cup } S^{\omega})$ where
$S^{\omega} = \{\gamma \in \Pi \, : \, \omega = \gamma \text{ or }
\omega - \gamma \in \Delta \}$. One may write $\mathfrak{r_{0}}=
\underset{\alpha \in \langle \Psi \rangle}{\sum
}\mathfrak{g}_{\alpha} + \underset{\alpha \in \langle \Psi
\rangle}{\sum }\mathbb{C}[e_{\alpha}, e_{-\alpha}]$ where $\langle
\psi \rangle $ is the set of roots which are linear combination of
elements of $\psi$. We prove in theorem \ref{theorem3.3} that the
appui subspace $V_{\tau}$ has the following form :

\hspace{1.5cm} $V_{\tau} = V_{\mathfrak{m}}$ \; \text{if } \;
$P_{\Psi}^{-}= \varnothing$

\hspace{1.5cm} $V_{\tau} = \mathfrak{g}$ \; \text{if }\;
$P_{\Psi}^{-}\neq \varnothing$ \vspace{0.5cm}

where $P_{\Psi}^{-}=\underset{\alpha \in \langle \Psi \rangle}{\sum
}\; \underset{\beta \in \Delta_{+} \backslash \mathfrak{R}_{2}}{\sum
} [\mathfrak{g}_{\alpha}, \mathfrak{g}_{-\beta}]$ and
$V_{\mathfrak{m}}$ is the appui subspace associated to
$\mathfrak{m}$.

\section{Standard subalgebras of affine Kac-Moody algebra}

Let $\mathfrak{g}$ be a simple Lie algebra of rank $p$. Let $A$ be
its extended Cartan matrix of affine type and $\mathfrak{g}(A)$ be
the affine Kac-Moody algebra associated to $A$. We consider the
element $\delta$ defined by $$ \delta =
\underset{i=0}{\overset{i=p}{\sum}}a_{i}\alpha_{i} $$

Let $\Delta^{aff}= \Delta^{re} \cup \Delta^{im}$ be a roots system
of the couple $(\mathfrak{g}(A), \mathfrak{h}+ \mathbb{C} \mathrm{K}
+ \mathbb{C} \texttt{d})$ where $\Delta^{re}=\{\alpha + j \delta \;
: \; \alpha \in \Delta \; \text{ and } \; j \in \mathbb{Z} \}$ is
the set of real roots and $ \Delta^{im}=\{j \delta \; : \; j \in
\mathbb{Z} \}$ is the set of imaginary roots.

We denote by $[.,.]$ the bracket on $\mathfrak{g}$. The bracket
$[.,.]_{aff}$ on $\mathfrak{g}(A)$ is defined as follows : for all
$(x,y) \in \mathfrak{g} \times \mathfrak{g} \; , \; (\lambda, \mu,
\lambda_{1}, \mu_{1}) \in \mathbb{C}^{4}\; , \; (n,m) \in
\mathbb{Z}^{2}$ :
$$ [t^{n} \otimes x + \lambda \mathrm{K} + \mu \texttt{d}, t^{m} \otimes y + \lambda_{1} \mathrm{K}
+ \mu_{1} \texttt{d}]_{aff} = t^{n+m} \otimes [x,y] + \mu m t^{m}
\otimes y - \mu_{1}n t^{n} \otimes x + n
\delta_{n,-m}(x|y)\mathrm{K}$$

In the particular case where $m$ and $n$ are two natural numbers, we
have : 
$$[t^{n} \otimes x + \lambda \mathrm{K} + \mu \texttt{d}, t^{m}
\otimes y + \lambda_{1} \mathrm{K} + \mu_{1} \texttt{d}]_{aff} =
t^{n+m} \otimes [x,y] + \mu m t^{m} \otimes y - \mu_{1}n t^{n}
\otimes x$$

\begin{definition} We introduce the following definitions :

\textit{$\mathrm{1.}$ A subalgebra of $\mathfrak{g}(A)$ is called
parabolic if it contains a Borel subalgebra of $\mathfrak{g}(A)$}.
\vspace{0.25cm}

\textit{$\mathrm{2.}$ A subalgebra of $\mathfrak{g}(A)$ is called
standard if its normalizer is a parabolic subalgebra of
$\mathfrak{g}(A)$}.
\end{definition}

Let $\tau$ be a standard subalgebra of $\mathfrak{g}$ of normalizer
$\rho(\tau)$ and let $V_{\tau}=[\tau, \mathfrak{g}]$ be the
associated subspace.

\begin{proposition}
\label {proposition1} The appui subspace $V_{\tau}$ satisfies the
following relation : $[V_{\tau}, \rho(\tau)] \subset V_{\tau}.$
\end{proposition}

\begin{proof}
Using the Jacobi identity, we have $[V_{\tau}, \rho(\tau)]= [[\tau,
\mathfrak{g}], \rho(\tau)] = [\tau, [\mathfrak{g}, \rho(\tau)]] +
[[\tau, \rho(\tau)], \mathfrak{g}] \subset [\tau, \mathfrak{g}] +
[\tau, \mathfrak{g}] \subset [\tau, \mathfrak{g}]=V_{\tau}$.
\end{proof}
\vspace{0,25 cm}

The following main theorem gives a sufficient condition for the
existence of a class of standard subalgebras of $\mathfrak{g}(A)$.

\begin{theorem} \label {existencetheorem}
Let $\tau$ be a standard subalgebra of $\mathfrak{g}$ of normalizer
$\rho(\tau)$ and $V_{\tau}=[\tau, \mathfrak{g}]$ be the appui
subspace associated to $\tau$.

Then, the subalgebra $$  \overline{\tau} = t^{n} \otimes \tau +
t^{n+1} \otimes V_{\tau} + t^{n+2}\mathbb{C}[t] \otimes \mathfrak{g}
+ \mathbb{C}K \; \; (n \in \mathbb{N}^{*}) $$  is a standard
subalgebra of $\mathfrak{g}(A)$ with normalizer
$\rho(\overline{\tau})= \rho(\tau) + t \mathbb{C}[t] \otimes
\mathfrak{g} + \mathbb{C}K + \mathbb{C}d$.
\end{theorem}

\begin{proof}
We show that the parabolic subalgebra $\rho(\overline{\tau})$ is the
normalizer of $\overline{\tau}$.

In $\mathfrak{g}(A)$, the subalgebra $\overline{\tau}$ is an ideal
of parabolic subalgebra $\rho(\tau) + t \mathbb{C}[t] \otimes
\mathfrak{g}+ \mathbb{C} \mathrm{K} + \mathbb{C} \texttt{d}$.

Let $\rho_{1}$ be a parabolic subalgebra of $\mathfrak{g}$ such that
$\rho_{1} + t\mathbb{C}[t] \otimes \mathfrak{g}+ \mathbb{C}
\mathrm{K} + \mathbb{C} \texttt{d}$ is a normalizer of
$\overline{\tau}$. We have $[\overline{\tau}, \rho_{1} +
t\mathbb{C}[t] \otimes \mathfrak{g}+ \mathbb{C} \mathrm{K} +
\mathbb{C} \texttt{d}]_{aff} \subset \overline{\tau}$ then $[\tau,
\rho_{1}] \subset \tau$ and $\rho_{1} \subset \rho(\tau)$. Since
$\rho_{1}+t \mathbb{C}[t]\otimes \mathfrak{g}+ \mathbb{C} \mathrm{K}
+ \mathbb{C} \texttt{d}$ is a normalizer of $\overline{\tau}$ then,
we have $\rho_{1}+t \mathbb{C}[t]\otimes \mathfrak{g} + \mathbb{C}
\mathrm{K} + \mathbb{C} \texttt{d}$ included in $\rho(\tau)+t
\mathbb{C}[t]\otimes \mathfrak{g}+ \mathbb{C} \mathrm{K} +
\mathbb{C} \texttt{d}$ and $\rho(\tau)$ included in $\rho_{1}$.
Furthermore, $\rho_{1}=\rho(\tau)$ and $\rho(\tau) + t \mathbb{C}[t]
\otimes \mathfrak{g}+ \mathbb{C} \mathrm{K} + \mathbb{C} \texttt{d}$
is a normalizer of $\overline{\tau}$.
\end{proof}

\begin{corollary}
The subalgebra $$\tau + t \otimes \mathfrak{g} + t^{2}\mathbb{C}[t]
\otimes \mathfrak{g}+ \mathbb{C}K + \mathbb{C}d$$  is a standard
subalgebra of $\mathfrak{g}(A)$ of normalizer $\rho(\tau) +
t\mathbb{C}[t] \otimes \mathfrak{g} + \mathbb{C}K + \mathbb{C}d$.
\end{corollary}

\begin{corollary}
Let $V$ be a subspace of $\mathfrak{g}$ such that $V_{\tau} \subset
V$ and $[V,\rho(\tau)] \subset V$ then the subalgebra
$$ t^{n} \otimes \tau + t^{n+1} \otimes V + t^{n+2}\mathbb{C}[t]
\otimes \mathfrak{g} + \mathbb{C}K \; \; (n \in \mathbb{N}^{*})$$ is
standard subalgebra of $\mathfrak{g}(A)$.
\end{corollary}

\begin{remark}
We deduce from the theorem \ref {existencetheorem} that the
subalgebra \; $\overline{V}_{\tau} = t^{n+1} \otimes V_{\tau} +
t^{n+2}\mathbb{C}[t] \otimes \mathfrak{g} + \mathbb{C}K$ \; is also
standard subalgebra of $\mathfrak{g}(A)$.
\end{remark}

\begin{example}
Let $\mathfrak{g}$ be a simple Lie algebra of type $B_{4}$ and $\Pi=
\{\alpha_{1}, \alpha_{2}, \alpha_{3}, \alpha_{4}\}$ be a basis of
roots system $\Delta$ .
\newline We set $\mathfrak{R}_{1} =
\{\alpha_{2} + 2\alpha_{3}+ 2\alpha_{4}, \alpha_{1} + \alpha_{2} +
2\alpha_{3}+ 2\alpha_{4}, \alpha_{1} + 2 \alpha_{2} + 2\alpha_{3}+
2\alpha_{4}\}$. We consider $\tau = \underset{\alpha \in
\mathfrak{R}_{1}}{\sum }\mathfrak{g}_{\alpha}$ and $V_{\tau} =
\underset{\alpha \in \Delta_{+} \backslash \{\alpha_{4}\}}{\sum
}\mathfrak{g}_{\alpha} + \mathfrak{h} + \mathfrak{g}_{-\alpha_{1}} +
\mathfrak{g}_{-\alpha_{2}} + \mathfrak{g}_{-\alpha_{1}-\alpha_{2}}$
respectively the standard subalgebra of $\mathfrak{g}$ and the appui
subspace associated to $\tau$.
\smallskip

Then $\overline{\tau} = t \otimes \tau + t^{2} \otimes V +
t^{3}\mathbb{C}[t] \otimes \mathfrak{g}+\mathbb{C}\mathrm{K}$ is a
standard subalgebra of $\mathfrak{g}(A)$ with normalizer
$\rho(\overline{\tau})=(\underset{\alpha \in \Delta_{+}}{\sum
}\mathfrak{g}_{\alpha} + \mathfrak{h} + \mathfrak{g}_{-\alpha_{1}}+
\mathfrak{g}_{-\alpha_{2}}+
\mathfrak{g}_{-\alpha_{1}-\alpha_{2}}+\mathfrak{g}_{-\alpha_{4}}) +
t\mathbb{C}[t] \otimes
\mathfrak{g}+\mathbb{C}\mathrm{K}+\mathbb{C}\texttt{d}$.

The standard subalgebra $\overline{\tau}$ may be written using the
real and imaginary roots as follows :$$\overline{\tau} =
\underset{\alpha \in \mathfrak{R}_{1}}{\sum} \mathfrak{g}(A)_{\alpha
+ \delta} + \underset{\alpha \in \mathfrak{R}_{\tau}}{\sum}
\mathfrak{g}(A)_{\alpha + 2 \delta} + \mathfrak{g}(A)_{2\delta} +
\underset{n \geqslant 3}{\sum } \{ \underset{\alpha \in
\Delta}{\sum} \mathfrak{g}(A)_{\alpha + n\delta} +
\mathfrak{g}(A)_{n\delta} \}+ \mathbb{C} \mathrm{K}$$ with
$\mathfrak{R}_{\tau}= \Delta_{+} \backslash \{\alpha_{4}\} \cup
\{-\alpha_{1}, -\alpha_{2}, -\alpha_{1}-\alpha_{2}\}$,
$\mathfrak{g}(A)_{\alpha + n \delta} = t^{n} \otimes
\mathfrak{g}_{\alpha}$ and $\mathfrak{g}(A)_{n \delta}=t^{n} \otimes
\mathfrak{h}$
\end{example}

Theorem \ref{existencetheorem} proves that the existence of this
class depends on the description of the standard subalgebras $\tau$
and the appui subspaces $V_{\tau}$ of the Lie algebra
$\mathfrak{g}$.

In the next section, we characterize the appui subspaces $V_{\tau}$
when $\tau$ is the nilpotent standard subalgebra then when $\tau$ is
the standard subalgebra.

\section{Appui subspaces}

Let $\mathfrak{g}$ be a complex simple Lie algebra, $\mathfrak{h}$
be its Cartan subalgebra and $\Delta$ be the roots system of the
couple $(\mathfrak{g}, \mathfrak{h})$. Let $\tau$ be the standard
subalgebra of $\mathfrak{g}$. The aim of this section is to give an
explicit formulas for the appui subspaces $V_{\tau}=[\tau,
\mathfrak{g}]$.

\subsection{Case of nilpotent standard subalgebra}

\textbf{Main notations.}
\newline $1$. For any positive root $\beta$, we denote by $S^{\beta}$ the set of simple roots
$\gamma$ of $\Pi$ such that $\beta = \gamma$ or $\beta - \gamma \in
\Delta$.  The set $S^{\beta}$ is called set of extremal roots of
$\beta$.
\newline $2$. For any root $\alpha=\underset{\gamma \in \Pi}{\sum} \alpha_{\gamma} \gamma$, we denote
by $C_{\alpha}$ the set of the simple roots $\gamma$ in $\Pi$ such
that $\alpha_{\gamma} \neq 0$.
\newline $3$. Let $B$ be the subset of the simple roots $\Pi$. We denote by $\langle B \rangle $ the set of roots
which are linear combination of elements of $B$ and by $\langle B
\rangle ^{+}$ (resp. $\langle B \rangle ^{-}$) the positive roots
(resp. the negative roots) of $\langle B \rangle$.
\medskip

Let $\tau$ be a nilpotent standard subalgebra of $\mathfrak{g}$
associated to the subsystem $\mathfrak{R}$ of positive roots two by
two non comparable. We set
\begin{equation} \label {R1} \mathfrak{R_{1}}= \{\, \beta \in
\Delta_{+} : \quad \omega \leqslant \beta \quad \text{for a certain
}\omega \in \mathfrak{R}\}\end{equation} The subalgebra $\tau$ has
the form $\tau = \underset{\beta \in \mathfrak{R}_{1}}{\sum
}{\mathfrak{g}_{\beta}}$.

We set $S_{2}=\underset{\omega \in \mathfrak{R}}{\cup }S^{\omega }$,
$\Delta _{1}= \langle \Pi \backslash S_{2} \rangle $ and $\Delta
_{2}=\Delta \backslash \Delta _{1}$. The normalizer of $\tau$ is
defined by $S_{2}$ and has the following form $\mathfrak{\rho
(\tau)}=\underset{\alpha \in \Delta _{+}}{\sum }\mathfrak{g}_{\alpha
}+\mathfrak{h} +\underset{\alpha \in \Delta _{1}^{+}}{\sum
}\mathfrak{g}_{-\alpha }$.

We can write $\mathfrak{g} = \rho(\tau) + \underset{\nu \in \Delta
_{2}^{+}}{\sum }\mathfrak{g}_{-\nu}$. Then, we have
\begin{equation}V_{\tau} = [\tau, \mathfrak{g}] \subset \tau +
\underset{\beta \in \mathfrak{R}_{1}}{\sum } \;\, \underset{\nu \in
\Delta_{2}^{+}}{\sum } [\mathfrak{g}_{\beta}, \mathfrak{g}_{-\nu}]
\label {equation1}.\end{equation}

\begin{lemma} \label{lemma2} Let $\alpha$ be a positive root.
Then $\mathfrak{g}_{\alpha}$ is included in $V_{\tau}$ if and only
if there exists $\beta \in \mathfrak{R}_{1}$ such that $C_{\alpha}
\cap S^{\beta} \neq \varnothing$.
\end{lemma}

\begin{proof}
($\Longrightarrow$) Let $\alpha \in \Delta_{+}$ such that
$\mathfrak{g}_{\alpha} \subset V_{\tau}$. There are two cases  :
\newline
First case : $\alpha \in \mathfrak{R}_{1}$ then there exists $\beta
\in \mathfrak{R}$ such that $\beta \leqslant \alpha$. This implies
that $C_{\alpha} \cap S^{\beta} \neq \varnothing$.
\newline
Second case : $\alpha \notin \mathfrak{R}_{1}$, using the formula
\ref {equation1}, there exists a root $\beta \in \mathfrak{R}_{1}$
such that $\beta - \alpha$ is a root of $\Delta_{2}^{+}$. One we can
write $\alpha$ in the form $\alpha = \alpha_{1}+ \cdots +
\alpha_{k}$ with $\alpha _{1}, \cdots, \alpha_{k}$ simple roots and
each partial sum $\alpha _{1} + \cdots + \alpha_{l}$ is a root, for
$1\leqslant l\leqslant k$. We want to prove by induction on $k$ that
there exists $i_{0} \in \llbracket{1,k} \rrbracket$ such that
$\alpha_{i_{0}}$ is in $S^{\beta}$.

If $k=1$. We have $\beta - \alpha_{1} \in \Delta$. Then $\alpha_{1}
\in S^{\beta}$. Suppose that $k
> 1$, then we have $ -\alpha _{k-1}^{'} - \alpha_{k} + \beta$ is a root where $\alpha
_{k-1}^{\prime} = \alpha _{1} + \cdots + \alpha_{k-1}$. By applying
Jacobi identity to the roots $-\alpha _{k-1}^{\prime}, -\alpha _{k}$
and $\beta$, we have $\beta - \alpha _{k-1}^{\prime}$ or $\beta -
\alpha_{k}$ is a root.
\newline If $\beta - \alpha_{k-1}^{\prime} \in \Delta$, by the assumption of recurrence, there
exists $j_{0} \in \llbracket{1,k-1} \rrbracket$ such that
$\alpha_{j_{0}} \in S^{\beta}$ and $C_{\alpha} \cap S^{\beta } \neq
\varnothing$.\newline If $\beta - \alpha_{k} \in \Delta$, we have
$\alpha_{k} \in S^{\beta}$ and $C_{\alpha} \cap S^{\beta } \neq
\varnothing$.\vspace{0,25cm}

($\Longleftarrow$) Let $\alpha$ be a root of $\Delta_{+}$ such that
there exists $\beta$ in $\mathfrak{R}_{1}$ with $C_{\alpha} \cap
S^{\beta} \neq \varnothing$. We set $\alpha = \alpha_{1} + \cdots +
\alpha_{n}$ such that for all $i \in \llbracket{1,n} \rrbracket$,
$\alpha_{i} \in \Pi$ and for $r \in \llbracket {1, n} \rrbracket$
each partial sum $\alpha_{1}+ \cdots + \alpha_{r} $ is a root.

Let $J$ be a subset of $\llbracket 1, n \rrbracket$ such that
$\alpha_{j} \in C_{\alpha} \cap S^{\beta}$ for any $j \in J$. It is
clear that the subset $J$ is not empty. Let $q$ be an element of
$J$. We have $\beta - \alpha_{q}$ is a root and
$\mathfrak{g}_{\alpha_{q}} = [\mathfrak{g}_{\beta},
\mathfrak{g}_{-\beta + \alpha_{q}}] \subset [\tau,
\mathfrak{g}]=V_{\tau}$.

Using the proposition \ref {proposition1}, we may deduce that
$$\mathfrak{g}_{\alpha}= \mathfrak{g}_{\alpha_{1}+ \cdots + \alpha_{n}} =
[[\mathfrak{g}_{\alpha_{q}}, \mathfrak{g}_{\alpha_{1} + \cdots +
\alpha_{q-1}}], \mathfrak{g}_{\alpha_{q+1}}, \cdots ,
\mathfrak{g}_{\alpha_{n}}] \subset [V_{\tau}, \rho(\tau)] \subset
V_{\tau}$$
\end{proof}

\begin{lemma}\label{lemma3} Let $\alpha$ be a positive root.
Then $\mathfrak{g}_{-\alpha}$ is included in $V_{\tau}$ if and only
if there exists $\beta \in \mathfrak{R}_{1}$ such that $\alpha +
\beta \in \Delta_{+}$.
\end{lemma}

\begin{proof}
($\Longrightarrow$) Let $\alpha \in \Delta_{+}$ such that
$\mathfrak{g}_{-\alpha} \subset V_{\tau}$. Then, by relation
\ref{equation1}, there exists a root $\beta \in \mathfrak{R}_{1}$
such that $\alpha + \beta$ is a root of $\Delta_{2}^{+}$.

($\Longleftarrow$) Let $\beta$ be a root in $\mathfrak{R}_{1}$ such
that $\alpha + \beta$ is a root. Then we have
$\mathfrak{g}_{-\alpha} = [\mathfrak{g}_{\beta},
\mathfrak{g}_{-\beta -\alpha}] \subset [\tau,
\mathfrak{g}]=V_{\tau}$.
\end{proof}
\medskip

We consider the set $\mathfrak{R}_{2}$ of $\Delta$ satisfying the
following property :
\begin{equation} \text{If } \alpha \in
\Delta_{+}, \beta \in \mathfrak{R}_{1}, \text{ and } \alpha + \beta
\in \Delta \text{ then } \alpha \in \mathfrak{R}_{2} \label {R2}
\end{equation}
And the set $\mathfrak{R}_{3}$ of $\Pi$ satisfying the following
property : \begin{equation} \text{If } \alpha \in \Pi, \beta \in
\mathfrak{R}_{1}, \text{ and } \beta - \alpha \notin \Delta \text{
then } \alpha \in \mathfrak{R}_{3} \label {R3} \end{equation} In
other words $\mathfrak{R}_{3} = \Pi \backslash S_{2}$

The next theorem characterizes the appui subspace $V_{\tau}$ when
$\tau$ is a nilpotent standard subalgebra of $\mathfrak{g}$.

\begin{theorem}
\label {theorem3.1} Let $\mathfrak{g}$ be a complex simple Lie
algebra and $\tau$ be a nilpotent standard subalgebra associated to
a subsystem $\mathfrak{R}$ of positive roots two by two non
comparable. Then, the appui subspace associated to $\tau$ has the
form
$$V_{\tau} = \underset{\alpha \in \Delta
_{+}\backslash <\mathfrak{R}_{3}>^{+}}{\sum }\mathfrak{g}_{\alpha }
+ \underset{\alpha \in \mathfrak{R}_{1}\cup \mathfrak{R}_{2}}{\sum}
\mathbb{C}\mathbb{[}e_{\alpha },e_{-\alpha }]+\underset{\alpha \in
\mathfrak{R}_{2}}{\sum }\mathfrak{g}_{-\alpha }$$
\end{theorem}

\begin{proof}
By applying the lemma \ref{lemma2} and \ref{lemma3}, it is enough to
prove that $\underset{\alpha \in \mathfrak{R}_{1}\cup
\mathfrak{R}_{2}}{\sum} \mathbb{C}\mathbb{[}e_{\alpha },e_{-\alpha
}]$ is included in $V_{\tau}$. \newline If $\alpha \in
\mathfrak{R_{1}}$ then $\mathbb{C}[e_{\alpha}, e_{-\alpha}] \subset
[\tau, \mathfrak{g}]=V_{\tau}$. If $\alpha \in \mathfrak{R}_{2}$
then, we have $\mathfrak{g}_{-\alpha}$ is included in $V_{\tau}$ and
by applying the proposition \ref{proposition1}, we have
$\mathbb{C}[e_{-\alpha}, e_{\alpha}] \subset [V_{\tau}, \rho(\tau)]
\subset V_{\tau}$.
\end{proof}

\begin{remark}
\label {remark3.3} We denote by $\theta$ the highest root of
$\Delta$. If the $\mathfrak{R}$is included in basis $\Pi$, the
nilpotent standard subalgebra $\tau$ is called complete standard
subalgebra of $\mathfrak{g}$ and it is easy to prove that the appui
subspace $V_{\tau}$ has the form
$$V_{\tau}= \underset{\alpha \in \Delta_{+}}{\sum
}\mathfrak{g}_{\alpha} + \mathfrak{h} + \underset{\alpha \in
\Delta_{+} \backslash \mathfrak{R}_{C}}{\sum
}\mathfrak{g}_{-\alpha}$$ with $\mathfrak{R}_{C} = \{\alpha \in
\Delta_{+} : \text{ for all root } \beta \in \mathfrak{R} \; , \;
\alpha_{\beta} = \theta_{\beta}\}$.
\end{remark}

\begin{remark}
We consider $m=\underset{\alpha \in
\Delta}{\sum}\mathfrak{g}_{\alpha}$ the nilpotent standard
subalgebra associated to the basis $\Pi$. It is the maximal
ad-nilpotent ideal of Borel subalgebra of $\mathfrak{g}$. Then, the
appui subspace is
$$V_{\tau}= \underset{\alpha \in \Delta_{+}}{\sum
}\mathfrak{g}_{\alpha} + \mathfrak{h} + \underset{\alpha \in
\Delta_{+} \backslash \{\theta\}}{\sum }\mathfrak{g}_{-\alpha}$$
\end{remark}

\subsection{General case }

In this section, we establish an explicit formula for the appui
subspaces $V_{\tau}$ when $\tau$ is a standard subalgebra of
$\mathfrak{g}$ not necessary nilpotent.

Let $\tau= \mathfrak{m} + \mathfrak{r}_{0}$ be a standard subalgebra
of $\mathfrak{g}$ of nilpotent radical where $\mathfrak{m}$ is its
nilpotent part and $\mathfrak{r}_{0}$ is its semisimple part. The
nilpotent part $\mathfrak{m}$ of $\tau$ is a nilpotent standard
subalgebra of $\mathfrak{g}$, then $\mathfrak{m}$ is associated to a
certain subset $\mathfrak{R}$ of positive roots two by two
noncomparable.

We recall that $S_{1} = \mathfrak{R} \cap \Pi$,
$S_{2}=\underset{\omega \in \mathfrak{R}}{\cup } S^{\omega}$ and
$\Delta_{1} = \langle \Pi \backslash S_{2} \rangle$. We set
$\Omega_{1} = \langle \Pi \backslash S_{1} \rangle$. We have
$\mathfrak{r}_{0}$ an ideal of the Levi reductive subalgebra
$\mathfrak{r_{1}} = \mathfrak{h} + \underset{\alpha \in
\Omega_{1}}{\sum }\mathfrak{g}_{\alpha}$ lying in the Levi reductive
subalgebra $\mathfrak{r_{2}} = \mathfrak{h} + \underset{\alpha \in
\Delta_{1}}{\sum }\mathfrak{g}_{\alpha}$. Since $\mathfrak{r}_{0}$
is semisimple algebra, then there exists $\Psi$ a common connected
subsystem of $\Pi \backslash S_{1}$ and $\Pi \backslash S_{2}$ such
that $\mathfrak{r_{0}}= \underset{\alpha \in \langle \Psi
\rangle}{\sum }\mathfrak{g}_{\alpha} + \underset{\alpha \in \langle
\Psi \rangle}{\sum }\mathbb{C}[e_{\alpha}, e_{-\alpha}]$.

\textit{Setting the subspace $P_{\Psi}^{-}=\underset{\alpha \in
\langle \Psi \rangle}{\sum } [\mathfrak{g}_{\alpha},
\mathfrak{n}_{2}^{-}]$ and the subspace
$P_{\Psi}^{+}=\underset{\alpha \in \langle \Psi \rangle}{\sum }
[\mathfrak{g}_{\alpha}, \mathfrak{n}_{2}^{+}]$} where
$\mathfrak{n}_{2}^{+}=\underset{\beta \in \Delta_{+} \backslash
\mathfrak{R}_{2}}{\sum }\mathfrak{g}_{\beta}$ and
$\mathfrak{n}_{2}^{-}=\underset{\beta \in \Delta_{+} \backslash
\mathfrak{R}_{2}}{\sum }\mathfrak{g}_{-\beta}$.

Using the properties of $\mathfrak{m}$ and $\mathfrak{r}_{0}$, one
proves that
$$V_{\tau} = V_{\mathfrak{m}} + P_{\Psi}^{-} +
[h_{\Psi},\mathfrak{n}_{2}^{-}]$$ with \hspace{0.05cm} $h_{\Psi} =
\underset{\alpha \in \langle \Psi \rangle}{\sum
}\mathbb{C}[e_{\alpha}, e_{-\alpha}]$ .

\begin{lemma}
\label {lemme3.2.1} If $\Psi$ is a common connected component of
$\Pi \backslash S_{1}$ and $\Pi \backslash S_{2}$ then
$P_{\Psi}^{+}$ is included in $\mathfrak{n_{2}^{+}}$.
\end{lemma}

\begin{proof}
Let $\alpha$ be a root such that $\mathfrak{g}_{\alpha} \subset
P_{\psi}^{+}$. Then, there exists $\gamma \in \langle \psi
\rangle^{+}$ and $\beta \in \Delta_{+} \backslash \mathfrak{R}_{2}$
such that $\alpha = \pm \gamma + \beta \in \Delta_{+}$. We want to
prove that $\mathfrak{g}_{\alpha} \subset \mathfrak{n}_{2}^{+}$.

We assume that $\alpha$ is a root of $\mathfrak{R}_{2}$. Then, there
exists $\omega \in \mathfrak{R}_{1}$ such that $\pm \gamma + \beta +
\omega \in \Delta$. Applying the Jacobi identity to $\pm \gamma$,
$\beta$ and $\omega$, we have two cases : \newline First case : if
$\pm \gamma + \omega \in \Delta$, since $\pm \gamma \in \langle \psi
\rangle \subset \Delta_{1}$ then $\pm \gamma +\omega \in
\mathfrak{R}_{1}$. We have $\underset{\in
\mathfrak{R}_{1}}{\underbrace{\pm \gamma + \omega}}+\beta$ then
$\beta \in \mathfrak{R}_{2}$. Therefor, we have a contradiction.
\newline Second case : if $\beta + \omega \in \Delta$ then $\beta
\in \mathfrak{R}_{2}$. Therefor, we have a contradiction.

Finally, $\alpha$ is a root of $\Delta_{+} \backslash
\mathfrak{R}_{2}$ and $\mathfrak{g}_{\alpha} \subset
\mathfrak{n}_{2}^{+}$.
\end{proof}

\begin{proposition}
\label {proposition3.2.2} Let $\alpha$ and $\beta$ be two roots of
$\Delta_{+} \backslash \mathfrak{R}_{2}$ such that $\alpha + \beta
\in \Delta_{+}$.
\newline
$1$. If $\mathfrak{g}_{\alpha+\beta } \subset P_{\Psi}^{+}$ then
$\mathfrak{g}_{\beta} \subset P_{\Psi}^{+}$ and
$\mathfrak{g}_{\beta} \subset P_{\Psi}^{+}$.
\newline
$2$. If $\mathfrak{g}_{\alpha} \subset P_{\psi}^{+}$ or
$\mathfrak{g}_{\beta} \subset P_{\psi}^{+}$ then
$\mathfrak{g}_{\alpha+\beta} \subset P_{\psi}^{+}$.
\end{proposition}

\begin{proof}
We have $\alpha + \beta \in \Delta_{+} \backslash \mathfrak{R}_{2}$
then, it is necessary that $\alpha \in \langle \Pi \backslash S_{1}
\rangle$ and $\beta \in \langle \Pi \backslash S_{1} \rangle$
because if no, we have $\alpha$ and $\beta$ are two roots of
$\mathfrak{R}_{1}$. So, by definition of the set $\mathfrak{R}_{2}$,
we have $\alpha$ and $\beta$ in $\mathfrak{R}_{2}$, contradiction.
\newline
$1$. We have $\mathfrak{g}_{\alpha + \beta} \subset P_{\psi}^{+}$
then, there exists $\gamma_{0} \in \langle \Psi \rangle^{+}$ such
that $\alpha + \beta \pm \gamma_{0} \in \Delta_{+} \backslash
\mathfrak{R}_{2}$. Applying the Jacobi identity to roots $\alpha$,
$\beta$ and $\pm \gamma_{0}$, we obtain two cases :
\newline First case : if $\alpha \pm \gamma_{0} \in
\Delta$ then $\alpha \pm \gamma \in \langle \psi \rangle$ because
$\psi$ is the connected subsystem of $\Pi \backslash S_{1}$.
Therefore, we have $[\mathfrak{g}_{-\alpha \mp \gamma_{0}},
\mathfrak{g}_{\beta + \alpha \pm \gamma_{0}}]= \mathfrak{g}_{\beta}
\subset P_{\psi}^{+}$ and $[\mathfrak{g}_{\mp \gamma_{0}},
\mathfrak{g}_{\alpha \pm \gamma_{0}}]= \mathfrak{g}_{\alpha} \subset
P_{\psi}^{+}$.
\newline Second case : if $\beta \pm \gamma_{0} \in \Delta$ then $\beta
\pm \gamma_{0} \in \langle \psi \rangle$ because $\psi$ is the
connected subsystem of $\Pi \backslash S_{1}$. Therefor, we have
$[\mathfrak{g}_{-\beta \mp \gamma_{0})}, \mathfrak{g}_{\alpha +
\beta \pm \gamma_{0}}]= \mathfrak{g}_{\alpha} \subset P_{\psi}^{+}$
and $[\mathfrak{g}_{\mp \gamma_{0}}, \mathfrak{g}_{\beta \pm
\gamma_{0}}]= \mathfrak{g}_{\beta} \subset P_{\psi}^{+}$.
\newline
\newline
$2$. Since $\mathfrak{g}_{\alpha} \subset P_{\psi}^{+}$ then, there
exists $\alpha_{0} \in \Delta \backslash \mathfrak{R}_{2}$ and
$\omega_{0} \in \langle \psi \rangle^{+}$ such that $\alpha =
\alpha_{0} \pm \omega_{0}$. Applying the Jacobi identity to $
\alpha_{0}$, $\pm \omega_{0}$ and $\beta$, we have two cases :
\newline First case : if $\alpha_{0} + \beta \in
\Delta_{+} \backslash \mathfrak{R}_{2}$ then $\mathfrak{g}_{\alpha +
\beta} \subset P_{\psi}^{+}$. \newline Second case : if $\beta \pm
\omega_{0} \in \Delta$ then $\beta \pm \omega_{0} \in \langle \psi
\rangle$ and $\mathfrak{g}_{\alpha + \beta} \subset P_{\psi}^{+}$.
\newline If $\mathfrak{g}_{\beta} \subset P_{\psi}^{+}$, we remplace
$\alpha$ by $\beta$ and with the same reason we prove that
$\mathfrak{g}_{\alpha+\beta} \subset P_{\psi}^{+}$.
\end{proof}

The next theorem gives the explicit formula of appui subspace
$V_{\tau}$ associated to standard subalgebra $\tau$ of
$\mathfrak{g}$.

\begin{theorem}
\label {theorem3.3} Let $\mathfrak{g}$ be a complex simple Lie
algebra and $\tau = \mathfrak{m} + \mathfrak{r_{0}}$ be a standard
subalgebra of $\mathfrak{g}$ of nilpotent radical, where
$\mathfrak{m}$ is a nilpotent standard subalgebra associated to a
subsystem $\mathfrak{R}$ of positive roots two by two non comparable
and $\mathfrak{r}_{0}$ is an ideal of $\mathfrak{r}_{1}$ contained
in $\mathfrak{r}_{2}$. \newline Then, we have : \newline \ 1.
$V_{\tau}=V_{\mathfrak{m}}$ \; if \; $P_{\Psi}^{-} = \varnothing$.
\newline \ 2. $V_{\tau}=\mathfrak{g}$ \; if \; $P_{\Psi}^{-} \neq
\varnothing $.
\end{theorem}

\begin{proof}

$1. $ If $P_{\Psi}^{-} = \varnothing$ then, $[h_{\Psi},
\mathfrak{n}^{-}_{2}] = \varnothing$. Therefore, $V_{\tau}=
[\tau,\mathfrak{g}] = V_{\mathfrak{m}}$ and we have the result.

$2.$ If $P_{\Psi}^{-} \neq \varnothing$. We have $V_{\tau} =
V_{\mathfrak{m}} + P_{\Psi}^{-} + [h_{\Psi},\mathfrak{n}_{2}^{-}]$
and $\mathfrak{g}=V_{\mathfrak{m}} + \mathfrak{n}_{2}^{-}$. So, it
is enough to show that $P_{\Psi}^{-} = \mathfrak{n_{2}^{-}}$. It is
the same to show that $P_{\Psi}^{+} = \mathfrak{n_{2}^{+}}$. By
using the lemma \ref {lemme3.2.1}, we have $P_{\Psi}^{+} \subset
\mathfrak{n_{2}^{+}}$, then it left to prove that
$\mathfrak{n_{2}^{+}} \subset P_{\Psi}^{+}$.\vspace{0,001cm}

We consider $\theta$ the highest root of $\Delta$. In the first, we
prove that $\mathfrak{g}_{\theta} \subset P_{\psi}^{+}$.
\newline We have $P_{\Psi}^{+} \neq \varnothing$ then, there
exists $\alpha$ a root of $\Delta_{+} \backslash \mathfrak{R_{2}}$
such that $\mathfrak{g}_{\alpha} \subset P_{\Psi}^{+}$. Since
$\alpha \leqslant \theta$, hence there exists $\theta_{1}, \cdots,
\theta_{k}$ a positive roots such that $\theta = \alpha + \theta_{1}
+ \cdots + \theta_{k}$ and $\alpha + \theta_{1} + \cdots +
\theta_{j} \in \Delta$ for each $j \in\llbracket{1, k} \rrbracket$.
By induction on $k$, we prove that $\mathfrak{g}_{\theta} \subset
P_{\Psi}^{+}$.
\newline If $k=1$, then $\theta = \alpha + \theta_{1}$ and by applying
the proposition \ref{proposition3.2.2}, we have
$\mathfrak{g}_{\theta} \subset P_{\psi}^{+}$. \newline If $k
\geqslant 1$. One may write $\theta
=\underset{\theta^{'}}{\underbrace{\alpha + \theta_{1} + \cdots +
\theta_{k-1}}} + \theta_{k}$. By induction
$\mathfrak{g}_{\theta^{'}} \subset P_{\psi}^{+}$ and by applying the
proposition \ref{proposition3.2.2}, we have $\mathfrak{g}_{\theta}
\subset P_{\psi}^{+}$.

Now, let $\beta$ be a root such that $\mathfrak{g}_{\beta} \subset
\mathfrak{n}_{2}^{+}$. We have $\beta \leqslant \theta$ then there
exists $\beta_{1}, \cdots, \beta_{k}$ a positive roots such that
$\theta = \beta + \beta_{1} + \cdots + \beta_{s}$ and $\beta +
\beta_{1} + \cdots + \beta_{j} \in \Delta$ for each $j
\in\llbracket{1, s} \rrbracket$. We have $\theta =
\underset{\beta^{'}}{\underbrace{\beta + \beta_{1} + \cdots +
\beta_{s-1}}} + \beta_{s}$, then by applying the proposition
\ref{proposition3.2.2}, we have $\mathfrak{g}_{\beta^{'}} \subset
P_{\psi}^{+}$. Then, we have $\mathfrak{g}_{\beta} \subset
P_{\psi}^{+}$. Therefore, we prove that $\mathfrak{n}_{2}^{+}
\subset P_{\psi}^{+}$ and this implies that $\mathfrak{n}_{2}^{+} =
P_{\psi}^{+}$ and $\mathfrak{n}_{2}^{-} = P_{\psi}^{-}$.

Finally, $[\tau, \mathfrak{g}]= V_{\mathfrak{m}} +
\mathfrak{n_{2}^{-}} = \mathfrak{g}$. So the second result of the
theorem is proved.
\end{proof}

\begin{example}
Let $\mathfrak{g}$ be a simple Lie algebra of type $F_{4}$ and $\Pi
= \{\alpha_{1}, \alpha_{2}, \alpha_{3}, \alpha_{4}\}$ be a basis of
the roots system of $F_{4}$. We consider the nilpotent standard
subalgebra $\mathfrak{m}$ associated to the subsystem $\mathfrak{R}
= \{\alpha_{3}\}$. We have $\Pi \backslash S_{2} = \Pi \backslash
S_{1} = \{\alpha_{1}, \alpha_{2}, \alpha_{4}\}$ and $\Delta_{+}
\backslash \mathfrak{R_{2}} = \{\alpha_{1} + 2\alpha_{2} +
4\alpha_{3} + 2\alpha_{4}, \alpha_{1} + 3\alpha_{2} + 4\alpha_{3} +
2\alpha_{4}, 2\alpha_{1} + 3\alpha_{2} + 4\alpha_{3} +
2\alpha_{4}\}$.
\begin{enumerate}
\item If we consider a standard subalgebra defined by $\mathfrak{m}$ and an ideal
$\mathfrak{r_{0}}$ associated to common connected component $\Psi =
\{\alpha_{4}\}$. We have $P_{\Psi}^{-} = \varnothing$, in this case
$V_{\tau} = V_{\mathfrak{m}}$.

\item If we consider a standard subalgebra defined by $\mathfrak{m}$ and an ideal
$\mathfrak{r_{0}}$ associated to common connected component $\Psi =
\{\alpha_{1}, \alpha_{2}\}$. We have $P_{\Psi}^{-} \neq
\varnothing$, more precisely, $P_{\Psi}^{-} = \mathfrak{n_{2}^{-}}$.
In this case $V_{\tau} = \mathfrak{g}$.
\end{enumerate}
\end{example}

\section{A class of standard subalgebra }

In this section, we prove that any standard subalgebra of
$\mathfrak{g}(A)$, under certain supplementary condition, is of the
form  given by theorem \ref{existencetheorem}.

Let $\mathfrak{m}$ be a nilpotent standard subalgebra of
$\mathfrak{g}$ associated to subsystem $\mathfrak{R} \subset
\Delta_{+}$ of roots two by two non comparable.

\begin{proposition}
\label {proposition3.4} Let $V_{\mathfrak{m}}$ be the appui subspace
associated to $\mathfrak{m}$ then $[V_{\mathfrak{m}},\mathfrak{g}] =
\mathfrak{g}$.
\end{proposition}

\begin{proof}
We consider $\mathfrak{n}$ the nilradical part of normalizer of
$\mathfrak{m}$. The subalgebra $\mathfrak{n}$ is complete standard
subalgebra associated to subsystem $S_{2}=\underset{\omega \in
\mathfrak{R}}{\cup } S^{\omega}$. We have $V_{\mathfrak{n}} =
[\mathfrak{n}, \mathfrak{g}] \subset [V_{\mathfrak{m}},
\mathfrak{g}]$. Since $V_{\mathfrak{n}} = \underset{\alpha \in
\Delta_{+}}{\sum }\mathfrak{g}_{\alpha} + \mathfrak{h} +
\underset{\alpha \in \Delta_{+ } \backslash
\mathfrak{R_{2}^{\mathfrak{n}}}}{\sum } \mathfrak{g}_{-\alpha}$ with
$\mathfrak{R}_{2}^{\mathfrak{n}} = \{\alpha \in \Delta_{+} :
\alpha_{\gamma} = \theta_{\gamma} \quad \text{for all } \gamma \in
S_{2}\}$ (remark \ref{remark3.3}) then, it is enough to prove that
$\mathfrak{g}_{-\alpha}$ is included in $[V_{\mathfrak{m}},
\mathfrak{g}]$, for any root $\alpha \in
\mathfrak{R}_{2}^{\mathfrak{n}}$. \newline Now, let $\beta$ be the
root of $\mathfrak{R}_{2}^{\mathfrak{n}}$ then for all $\gamma \in
S^{\beta}$ : $ht_{\gamma}(\alpha) \geqslant ht_{\gamma}(\beta)$
where $ht_{\gamma}(\alpha)$ and $ht_{\gamma}(\beta)$ is the
multiplicity of $\gamma$ in the decomposition of simple roots of
$\alpha$ and of $\beta$. Therefore, we have $\beta \leqslant \alpha$
and this implies that $\alpha \in \mathfrak{R_{1}}$. We deduce that
$\mathfrak{g}_{-\alpha} = [[e_{\alpha}, e_{-\alpha}],
\mathfrak{g}_{-\alpha}] \subset [V_{\mathfrak{m}}, \mathfrak{g}]$.
This proves that $V = \mathfrak{g}$.
\end{proof}

\begin{remark} \label{generalcase}
In proposition \ref{proposition3.4}, if we take $\tau$ a standard
subalgebra not necessarily nilpotent, we have also
$[V_{\tau},\mathfrak{g}] = \mathfrak{g}$.
\end{remark}

Let $n$ be a natural number. Let  $\{I_{n+j}\}_{j \geqslant 0}$ be a
family of vectors subspaces of $\mathfrak{g}$ such that the subspace
$\overline{\tau} = \underset{j \geqslant 0}{\sum }t^{n+j} \otimes
I_{n+j} + \mathbb{C}\mathrm{K}$ is a standard subalgebra of
$\mathfrak{g}(A)$ of normalizer $\overline{\rho } = \rho +
t\mathbb{C}[t] \otimes \mathfrak{g}+ \mathbb{C} \mathrm{K} +
\mathbb{C} \texttt{d}$ with $\rho$ is the parabolic subalgebra of
$\mathfrak{g}$. \vspace{0.25cm}

We have $[\overline{\tau}, \overline{\rho}] = t^{n} \otimes [I_{n},
\rho] + t^{n+1} \otimes ([I_{n}, \mathfrak{g}] + [I_{n+1}, \rho]) +
t^{n+2} \otimes ([I_{n}, \mathfrak{g}] + [I_{n+1}, \mathfrak{g}] +
[I_{n+2}, \rho]) + t^{n+3} \otimes (\cdots + [I_{n+2}, \mathfrak{g}]
+ \cdots) + \cdots \subset \overline{\tau}$.

\begin{lemma} \label{lemmarelation}
We have the following relations :
\begin{enumerate}
\item $[I_{n}, \rho] \subset I_{n}$.
\item $[I_{n}, \mathfrak{g}] \subset I_{n+1}$ and $[I_{n+1}, \rho] \subset I_{n+1}$.
 \item $[I_{n+1},
\mathfrak{g}] \subset I_{n+2}$ and $[I_{n+2}, \rho] \subset
I_{n+2}$.
\item For all $j \geqslant 3$ : $[I_{n+2},
\mathfrak{g}] \subset I_{n+j}$.
\end{enumerate}
\end{lemma}

\begin{proposition}
\label{importantpropo} If $I_{n}$ is a Lie subalgebra of
$\mathfrak{g}$ then $I_{n}$ is included in $\rho$.
\end{proposition}

\begin{proof}
We can write $\rho = \underset{\alpha \in \Delta_{+}}{\sum
}\mathfrak{g}_{\alpha} + \mathfrak{h} + \underset{\alpha \in \langle
\Pi \backslash T \rangle^{+}}{\sum }\mathfrak{g}_{-\alpha}$, for a
certain subset $T$ of $\Pi$.
\newline We suppose that $I_{n}$ is not in $\rho$ and
we set $T_{1} = \{ \beta \in T , \; \mathfrak{g}_{-\beta} \subset
I_{n} \}$. We consider $\rho_{1} = \underset{\alpha \in
\Delta_{+}}{\sum }\mathfrak{g}_{\alpha} + \mathfrak{h} +
\underset{\alpha \in \langle \Pi \backslash (T \backslash T_{1})
\rangle^{+}}{\sum }\mathfrak{g}_{-\alpha}$ the parabolic subalgebra
of $\mathfrak{g}$ associated to $(T \backslash T_{1})$. By
definition $\rho$ is included in $\rho_{1}$.

We want to prove that $\overline{\tau}$ is an ideal of $\rho_{1} + t
\mathbb{C}[t] \otimes \mathfrak{g}+ \mathbb{C} \mathrm{K} +
\mathbb{C} \texttt{d}$. \newline First, we prove that $I_{n}$ is an
ideal of $\rho_{1}$. Let $\alpha$ and $\beta$ be two roots such that
$\mathfrak{g}_{\alpha} \subset I_{n}$, $\mathfrak{g}_{\beta} \subset
\rho_{1}$ and $\alpha + \beta$ is the root. \newline If $\beta$ is a
root of $\Delta_{+} \cup \langle \Pi \backslash T \rangle^{-}$ then
$\mathfrak{g}_{\beta}$ is included in $\rho$ and
$[\mathfrak{g}_{\alpha}, \mathfrak{g}_{\beta}] \subset [I_{n}, \rho]
\subset I_{n}$. \newline If $\beta$ is a root of $\langle \Pi
\backslash (T \backslash T_{1}) \rangle^{-}$ with $C_{\beta} \cap
T_{1} \neq \varnothing$. There exists $\beta_{1}, \cdots, \beta_{r}$
a simple roots such that $\beta = -\beta_{1} - \cdots - \beta_{r}$
and for all $j \in \llbracket{1, r} \rrbracket$, $\beta_{1} + \cdots
+ \beta_{j} \in \Delta$. Let $i$ be the smallest index such that
$\beta_{i} \in T_{1}$. We have $\beta_{i-1}^{'}= \beta_{1} + \cdots
+ \beta_{i} \in \langle \Pi \backslash T\rangle$ and
$\mathfrak{g}_{-\beta^{'}_{i}} = [
\mathfrak{g}_{-\beta_{i}},\mathfrak{g}_{-\beta^{'}_{i-1}}] \subset
[I_{n},\rho] \subset I_{n}$. Now taking the simple root
$\beta_{i+1}$, we have two cases :

If $\beta_{i+1} \in T_{1}$, we have $\mathfrak{g}_{-\beta_{i+1}^{'}}
= [\mathfrak{g}_{-\beta_{i}^{'}}, \mathfrak{g}_{-\beta_{i+1}}]
\subset [I_{n}, I_{n}] \subset I_{n}$ (because $I_{n}$ is
subalgebra).

If $\beta_{i+1} \in \Pi \backslash T$, we have
$\mathfrak{g}_{-\beta_{i+1}^{'}} = [
\mathfrak{g}_{-\beta_{i+1}},\mathfrak{g}_{-\beta_{i}^{'}}] \subset
[\rho, I_{n}] \subset I_{n}$. \newline With similar argument, we
show that $\mathfrak{g}_{\beta} \subset I_{n}$.

Therefore, $\mathfrak{g}_{\alpha + \beta} = [\mathfrak{g}_{\alpha},
\mathfrak{g}_{\beta}] \subset [I_{n}, I_{n}] \subset I_{n} $. Then,
$I_{n}$ is the standard subalgebra of $\mathfrak{g}$ of normalizer
$\rho_{1}$.

Now, we have $V_{I_{\mathfrak{n}}}=[I_{n}, \mathfrak{g}] \subset
I_{n+1} $ and $[I_{n+1}, \rho_{1}] \subset [I_{n+1}, \rho] +
[I_{n+1}, I_{n}] \subset I_{n+1}$.

Since $I_{n}$ is a standard subalgebra then, by applying the
proposition \ref{proposition3.4} an the relation $3$ of the lamma
\ref{lemmarelation}, we have $\mathfrak{g} =[V_{I_{n}},
\mathfrak{g}] \subset [I_{n+1},\mathfrak{g}] \subset I_{n+2}$. Then
$I_{n+2}=\mathfrak{g}$.
\newline Since $\mathfrak{g}$ is the simple Lie
algebra, for $j \geqslant 3$, we have $[I_{n+2},
\mathfrak{g}]=[\mathfrak{g},\mathfrak{g}]= \mathfrak{g} \subset
I_{n+j}$. This implies that $I_{n+j}=\mathfrak{g}$.

Finally, we have proved  that $\overline{\tau}$ is an ideal of the
parabolic subalgebra $\rho_{1} + t \mathbb{C}[t] \otimes
\mathfrak{g}+\mathbb{C}\mathrm{K}+\mathbb{C}\texttt{d}$ and this
subalgebra contains the normalizer $\rho + t \mathbb{C}[t] \otimes
\mathfrak{g}+\mathbb{C}\mathrm{K}+\mathbb{C}\texttt{d}$ of
$\overline{\tau}$, contradiction. Hence $I_{n}$ is included in
$\rho$.
\end{proof}

\begin{remark}
The proposition \ref{importantpropo} proves that $I_{n}$ is the
standard subalgebra of normalizer $\rho$.
\end{remark}

The following theorem determines a class of standard subalgebra of
$\mathfrak{g}(A)$ under certain condition.

\begin{theorem}
Let $\overline{\tau} = \underset{j \geqslant 0}{\sum }t^{n+j}
\otimes I_{n+j}+ \mathbb{C}\mathrm{K}$ be a standard subalgebra of
$\mathfrak{g}(A)$ of normalizer $\overline{\rho}=\rho +
t\mathbb{C}[t] \otimes \mathfrak{g}+ \mathbb{C}\mathrm{K} +
\mathbb{C} \texttt{d}$ and  $I_{n} \neq 0$. \vspace{0,25cm}

If $I_{n}$ is a Lie subalgebra of $\mathfrak{g}$ then \;
$\overline{\tau} = t^{n} \otimes \tau + t^{n+1} \otimes V + t^{n+2}
\mathbb{C}[t]\otimes \mathfrak{g} + \mathbb{C}\mathrm{K}$ \quad
where $\tau = I_{n}$ is a standard subalgebra of $\mathfrak{g}$ of
normalizer $\rho$, $V$ is a subspace of $\mathfrak{g}$ contain the
subspace $V_{\tau}=[\tau, \mathfrak{g}]$ and $I_{n+j}=\mathfrak{g}$,
for all $j \geqslant 2$.
\end{theorem}

\begin{proof}
Let $n$ be a natural number. We set $I_{n}=\tau$. We write
$\overline{\tau} = t^{n} \otimes \tau + \underset{j \geqslant
1}{\sum }t^{n+j} \otimes I_{n+j}+ \mathbb{C}\mathrm{K}$.

By applying the proposition \ref{importantpropo}, we have $\tau$ is
included in $\rho$. Then, $\tau$ is the standard subalgebra of
$\mathfrak{g}$ of normalizer $\rho$.

We set $V = I_{n+1}$. The relation $2.$ in the lemma
\ref{lemmarelation} proves that $V_{\tau}$ is included in $V$ where
$V_{\tau}$ is the appui subspace associated to $\tau$.

We prove by induction on $j\geqslant 2 $ that $I_{n+j} =
\mathfrak{g}$. We use the relation $3$ in the lemma
\ref{lemmarelation}, then we have $[V_{\tau}, \mathfrak{g}] \subset
I_{n+2}$. By applying the proposition \ref{proposition3.4} and
remark \ref{generalcase}, we have $I_{n+2} = \mathfrak{g}$.

Since $\mathfrak{g}$ is the simple Lie algebra, then we have
$[\mathfrak{g}, \mathfrak{g}] = \mathfrak{g}$. Using the relation
$4$ in the lemma \ref{lemmarelation}, we have
$I_{n+j}=\mathfrak{g}$. Therefore $\overline{\tau} = t^{n} \otimes
\tau + t^{n+1} \otimes V + t^{n+2}\mathbb{C}[t] \otimes
\mathfrak{g}$. \vspace{0,10cm}
\end{proof}


\begin{corollary}
Let $\overline{\tau} = I_{0} + \underset{k \geqslant 1}{\sum}t^{k}
\otimes I_{k}+\mathbb{C}\mathrm{K}+\mathbb{C}\texttt{d}$ be the
standard subalgebra of $\mathfrak{g}(A)$ of normalizer
$\rho(\overline{\tau})=\rho + t \mathbb{C}[t] \otimes \mathfrak{g} +
\mathbb{C} \mathrm{K} + \mathbb{C} \texttt{d}$.

Then $\tau = I_{0}$ is a standard subalgebra of $\mathfrak{g}$ of
normalizer $\rho$ and $\overline{\tau}= \tau + t \mathbb{C}[t]
\otimes \mathfrak{g} + \mathbb{C} \mathrm{K} + \mathbb{C}
\texttt{d}$.
\end{corollary}

\begin{proof}
It is an immediate consequence of the previous theorem.
\end{proof}


\end{document}